\documentclass[a4paper, 10pt]{amsart}

\usepackage[latin1]{inputenc}

\usepackage{amssymb, amsmath}

\usepackage[english]{babel}

\usepackage[alphabetic]{amsrefs}

\numberwithin{equation}{section}

\newcommand{\Z}{{\mathbb Z}}

\newcommand{\KK}{{\mathbf K}}\newcommand{\QQ}{{\mathbf Q}}

\def\:{\colon}
\def\.{\cdot}
\def\o{\circ}
\def\<{\left\langle}
\def\>{\right\rangle}
\def\({\left(}
\def\){\right)}
\def\epsilon{\varepsilon}
\def\phi{\varphi}
\def\subset{\subseteq}

\def\leq{\leqslant}
\def\geq{\geqslant}

\def\lra{\longrightarrow}

\def\mapsto{\longmapsto}
\def\Mon{\mathbf{Mon}}
\def\A{\mathcal A}
\def\dd{\partial\partial}

\newtheorem{theorem}[equation]{Theorem}

\newtheorem{proposition}[equation]{Proposition}

\numberwithin{equation}{section}

\theoremstyle{remark}
\newtheorem{definition}[equation]{Definition}
\newtheorem{remark}[equation]{Remark}

\DeclareMathOperator{\Hom}{Hom} 
\DeclareMathOperator{\gr}{gr} \DeclareMathOperator{\Ext}{Ext}
\DeclareMathOperator{\Tor}{Tor} \DeclareMathOperator{\im}{im}
 \DeclareMathOperator{\Tot}{Tot}
\DeclareMathOperator{\coker}{coker}

\DeclareMathOperator{\Ch}{\mathbf{Ch}}
\DeclareMathOperator{\DC}{\mathbf{DC}}
\DeclareMathOperator{\twist}{\widetilde\otimes}

\title{Homology of powers of regular ideals}
\author{Samuel W\"uthrich}
\address{Mathematisches Institut, Universit\"at Bern,
Sidlerstrasse 5, 3012 Bern, Switzerland.}
\email{s.wuethrich@math-stat.unibe.ch} 

\date{27/08/03}

\begin{document}

\begin{abstract}
For a commutative ring $R$ with an ideal $I$, generated by a
finite regular sequence, we construct differential graded algebras
which provide $R$--free resolutions of $I^s$ and of $R/I^s$ for
$s\geq 1$ and which generalise the Koszul resolution. We derive
these from a certain multiplicative double complex $\KK$. By means
of a Cartan--Eilenberg spectral sequence we express $\Tor_*^R(R/I,
R/I^s)$ and $\Tor_*^R(R/I, I^s)$ in terms of exact sequences and
find that they are free as $R/I$--modules. Except for $R/I$, their
product structure turns out to be trivial; instead, we consider an
exterior product $\Tor_*^R(R/I, I^s)\otimes_R\Tor_*^R(R/I, I^t)\to
\Tor_*^R(R/I, I^{s+t})$. --- This paper is based on ideas by
Andrew Baker; it is written in view of applications to algebraic
topology.
\end{abstract}

\maketitle
\section*{Introduction}\label{introduction}
Let $R$ be a commutative ring with unit and let $I\lhd R$ be an
ideal generated by a finite regular sequence $r_1, \ldots, r_n\in
R$, i.e.\@ $r_1$ is a non-zero divisor of $R$ such that
$R/(r_1)\neq 0$, $r_2$ a non-zero divisor of $R/(r_1)$ such that
$R/(r_1, r_2)\neq 0$ and so on. The Koszul complex $K$, a
differential graded algebra, provides a canonical $R$--free
resolution of $R/I$. The aim of this paper is to construct
explicit $R$--free resolutions of $R/I^s$ and $I^s$ for $s>1$,
which generalise the Koszul resolution and which enable us to
compute $\Tor_*^R(R/I, R/I^s)$ and $\Tor_*^R(R/I, I^s)$.

We derive these from a certain fourth quadrant double complex
$\KK$ with $R$--free components, which we call the extended Koszul
complex. In formal analogy to the ordinary Koszul complex, $\KK$
is constructed as a tensor product of elementary double complexes
and carries a multiplicative structure. Its filtration by columns
turns out to be a filtration by ideals $F^s(\KK)\lhd\KK$, so that
the associated total complexes of $F^s(\KK)$ and
$\KK/s=\KK/F^s(\KK)$ are differential graded algebras. Via
suitable augmentations, they provide resolutions of $I^s$ and
$R/I^s$ respectively.

The column-wise filtration of $R/I\otimes_R\, F^s(\KK)$ and
$R/I\otimes_R\, \KK/s$ give rise to Cartan--Eilenberg spectral
sequences converging to $\Tor_*^R(R/I, R/I^s)$ and $\Tor_*^R(R/I,
I^s)$ respectively. We show that they collapse at $E^2$, express
the $\Tor$--groups by means of exact sequences and find that they
are free over $R/I$.

Except for $R/I$, their multiplicative structure, induced by the
$R$--algebra structure on $R/I^s$ and on $I^s$, turns out to be
trivial. The multiplications $I^s\otimes_R I^t \to I^{s+t}$
however induce non-trivial exterior products
\[\Tor_*^R(R/I, I^s)\otimes_R \Tor_*^R(R/I, I^t) \lra \Tor_*^R(R/I,
I^{s+t})\] which allow us to view $\bigoplus_{s\geq 0}
\Tor_*^R(R/I, I^s)$ as a bigraded algebra. It contains a copy of
the polynomial ring $R/I[x_1,\ldots, x_n]$ as a subalgebra, over
which the whole algebra is generated by the $R/I$--basis elements
of $\Tor_{k>0}^R(R/I,I)$.

This paper is based on ideas by Andrew Baker; he constructs in
\cite{B} resolutions of $R/I^s$ and computes $\Tor_*^R(R/I,
R/I^s)$. His approach to the resolutions, however, is conceptually
different from the one presented here. He constructs them
inductively, by pasting together copies of the Koszul complex.

The present paper originated in the course of the work on my PhD
thesis, in view of applications to algebraic topology. These will
be discussed in \cite{W}.

Note that Tate constructs in \cite{T} for \emph{any} $R$--algebra
of the form $R/M$, where $R$ is a commutative Noetherian ring with
unit and $M\lhd R$ an ideal, a differential graded algebra which
provides a free resolution. Applying the construction to $R/I^s$,
where $I$ is as above, does not give the resolution considered
here.

It is a pleasure to express my gratitude to my supervisor Alain
Jeanneret from the University of Berne and Andrew Baker from the
University of Glasgow for many interesting discussions and for
ongoing support. In particular, I would like to thank Andrew Baker
for inviting me to spend a year at the University of Glasgow, as
well as the Mathematical Institute in Berne for giving me the
opportunity to do so. I would like to thank the people from the
Mathematical Department in Glasgow in general, for offering an
inspiring and stimulating atmosphere during my stay.

\section{Multiplicative double complexes}\label{multiplicative}

We define the notion of ``double complex'' which is most suitable
for our purposes and describe a natural way of forming tensor
products; this allows us to consider multiplicative double
complexes. References for background material for this section are
\cite{M1}*{VII.}, \cite{M2}*{X.9} and \cite{We}.

Given an abelian category $\A$, we can form the category $\Ch(\A)$
of chain complexes in $\A$, whose morphisms are the chain maps. As
$\Ch(\A)$ in turn is abelian, in a canonical way, we can iterate
the construction.
\begin{definition}
A {\em double complex} is an object $C$ of $\DC(\A)=
\Ch(\Ch(\A))$.
\end{definition}

We display the components $C_p$ of a double complex $C$ as the
columns of a lattice $C_{p,q}$. By definition, the differential
$d_p\:C_p\to C_{p-1}$ is a chain map, so the squares in the
lattice commute. Note that one often means by a double complex a
lattice $C_{p,q}$ whose rows and columns are chain complexes, as
here, but whose squares {\em anti}commute.

We denote the components $C_{p,q} \to C_{p-1,q}$ of the
differential of $C$ by $d_{p,q}^h$ (the {\em horizontal
differentials}\/) and the components $C_{p,q} \to C_{p,q-1}$ of
the differentials of the columns $C_p$ of $C$ by $d_{p,q}^v$ (the
{\em vertical differentials}\/); when the indices are clear from
the context, we omit them. We also use the convention $C^p_q=
C_{-p,q}$.

By replacing the vertical differential $d^v_{p,q}$ by
$(d^v_{p,q})' = (-1)^p d_{p,q}^v$, we can pass from a double
complex $C$, as defined above, to one with anticommuting squares,
which we denote by $C^\o$.

The total complex $\Tot^\oplus(C^\o)$ of $C^\o$ has components
\[(\Tot^\oplus(C^\o))_k = \bigoplus_{p+q=k} C_{p,q}\]
and differential given by $d= d^h + (d^v)'$. Following MacLane, we
call this complex the {\em condensation} of $C$ and denote it by
$C^\bullet$.

Assume now that $\< \A, \otimes, e\>$ is an abelian, symmetric
monoidal category (where $e$ is the unit). We will refer to the
bifunctor $\otimes$ as the ``tensor product''.

Recall that $\Ch(\A)$ inherits a symmetric monoidal structure from
$\A$. Namely, for two chain complexes $(C,d)$ and $(D,d')$, the
tensor product $C\otimes D$ is defined by
\begin{equation}\label{tensorproduct}
 (C\otimes D)_k = \bigoplus_{p+q=k} C_p\otimes D_q,
\end{equation}
with differential given by
\[d_p \otimes 1 + (-1)^p\, 1\otimes d'_q\: C_p\otimes
D_q \lra (C\otimes D)_{p+q-1}.\] Embedding $\A$ in $\Ch(\A)$ in
the usual way, $e$ is a unit in $\Ch(\A)$ for $\otimes$. The
symmetry isomorphism $\Hat\tau: C\otimes D \lra D\otimes C$ is
defined as
\begin{equation}\label{symmetry}
 \Hat\tau(c\otimes d) = (-1)^{kl} \tau(c\otimes d)
\end{equation}
for $c\otimes d\in C_k\otimes D_l$, where $\tau: C_k\otimes
D_l\lra D_l\otimes C_k$ is the given symmetry isomorphism in $\A$.

Iterating this procedure, we get a symmetric monoidal structure on
$\DC(\A)$. So by definition, the components (columns) of the
tensor product $C\otimes D$ of two double complexes $C$ and $D$
consist of a direct sum of tensor products of columns of $C$ and
$D$. Note that $C\otimes D$ has the analogous property for the
rows. We can express this in a more conceptual way by introducing
the {\em transpose} $TC$ of a double complex $C$, defined as
$(TC)_{p,q}= C_{q,p}$, with differentials $(Td)^v_{p,q}=
d^h_{q,p}$ and $(Td)^h_{p,q}=d^v_{q,p}$. The statement then
amounts to the equation
\begin{equation}\label{rows}
TC\otimes TD = T(C\otimes D).
\end{equation}
Note that the symmetry isomorphism $\Hat{\Hat{\tau}}:C\otimes D\to
D\otimes C$ is given by \[\Hat{\Hat{\tau}}(c\otimes d) =
(-1)^{kp+lq}\tau(c\otimes d)\] for $c\otimes d\in C_{k,l}\otimes
D_{p,q}$.

\begin{definition}
A {\em multiplicative double complex} is a monoid in the monoidal
category $\< \DC(\A),\otimes,e \>$, i.e.\@ a double complex $C$
with a multiplication $\mu\: C\otimes C \to C$ and a unit $\eta\:
e\to C$ such that the associativity and the two unit diagrams
commute.
\end{definition}

Unravelling the definitions, the product $\mu$ of a multiplicative
double complex is defined by a collection of maps $C_{p,q}\otimes
C_{k,l}\to C_{p+k, q+l}$ such that both the components of the
horizontal and the vertical differential are derivations, in the
sense that we have, for $c\in C_{p,q}, \ d\in C_{k,l}$,
\begin{alignat}{2}\label{diffbigradedh}
d^h(\mu(c\otimes d)) &= \mu(d^h(c), d) + (-1)^p\, \mu(c, d^h(d))
\\ \label{diffbigradedv} d^v(\mu(c\otimes d))
&= \mu(d^v(c), d) + (-1)^q\, \mu(c, d^v(d)).
\end{alignat}

Note that the condensation $C^\bullet$ of a multiplicative double
complex is canonically a monoid in $\Ch(\A)$, as a consequence of
the natural isomorphism $(C\otimes D)^\bullet \cong
C^\bullet\otimes D^\bullet$, which holds for any two double
complexes $C$ and $D$, and the fact that condensation is a
functorial process.

Suppose that $(C, \mu, \eta)$ and $(D, \mu', \eta')$ are two
multiplicative double complexes. As in any symmetric monoidal
category, we can endow $C\otimes D$ with a multiplicative
structure in a canonical way, by defining the product as the
composition
\begin{equation}\label{monoid}
(C\otimes D)\otimes(C\otimes D) \xrightarrow{C\otimes
\Hat{\Hat{\tau}} \otimes D} C\otimes C\otimes D\otimes D
\xrightarrow{\mu\otimes \mu'} C\otimes D.
\end{equation}
Explicitly, the product is given by the collection of maps
\begin{align*}
(C_{p,q}\otimes D_{k,l}) \otimes (C_{p',q'}\otimes D_{k',l'}) &
\lra C_{p+p', q+q'}\otimes D_{k+k', l+l'}\\
(c\otimes d)\otimes(c'\otimes d') & \mapsto (-1)^{k p' + l q'}
\mu(c\otimes c') \otimes \mu'(d\otimes d').
\end{align*}
We have morphisms of monoids
\begin{align*}
 C\cong C\otimes e \xrightarrow{C\otimes\eta'} C\otimes D, \quad
D\cong e\otimes D \xrightarrow{\eta\otimes D} C\otimes D.
\end{align*}

We will also need the category $\dd\A$ of ($\Z$--)bigraded objects
in $\A$; by the construction of \eqref{tensorproduct}, $\dd\A$ is
monoidal. The forgetful functor $U\: \DC(\A)\to\dd\A$ is a strict
morphism of monoidal categories and therefore restricts to a
functor
\[ U\: \Mon_{\DC(\A)} \lra \Mon_{\dd\A}\]
between the categories of monoids. Whereas in $\DC(\A)$ we were
forced to introduce a sign when defining the symmetry isomorphisms
in \eqref{symmetry}, we wouldn't need to do so in $\dd\A$.
However, we want the restriction of $U$ to $\Mon_{\DC(\A)}$ to be
monoidal as well; therefore, we define symmetry isomorphisms in
$\dd\A$ as in \eqref{symmetry}. To stress that we use these
symmetries for the definition of the tensor product of monoids in
$\dd\A$, as in \eqref{monoid}, we denote it by $\twist$ and call
it the {\em twisted tensor product}.

In the next section, we apply these constructions to the symmetric
monoidal category $\< R\mathbf{-Mod},\otimes_R, R\>$ of
$R$--modules. The symmetry isomorphisms are given by the switch
maps $\tau: M\otimes N\to N\otimes M, \, \tau(m\otimes n)=
n\otimes m$. Monoids in $R\mathbf{-Mod}$ and $\Ch(R\mathbf{-Mod})$
are $R$--algebras and differential graded $R$--algebras
respectively. We will refer to monoids in $\dd(R\mathbf{-Mod})$ as
{\em bigraded $R$--algebras.}

\section{The extended Koszul complex and the resolutions}
\label{extendedkoszul}

Without further notice, the ground ring is from now on understood
to be $R$, so we omit the letter $R$ in $\otimes_R, \ \Hom_R,
\Lambda_R, \ \Tor^R$ or $\Ext_R$.

We briefly recall the definition of the Koszul complex $K$, in
order to fix notations and to point out the formal similarity to
the double complex $\KK$ we construct afterwards.

For a non-zero divisor $r\in R$, $K(r)$ is defined to be the
differential graded algebra $(\Lambda(e),d)$ with $|e|=1$ and
$d(e)=r$. The projection $R\to R/(r)$ defines an augmentation
$\epsilon_r\:K(r)\to R/(r)$, which exposes $K(r)$ as an $R$--free
resolution of $R/(r)$. For a regular sequence $r_1, \ldots, r_n$
generating an ideal $I$, the Koszul complex $K$ is the
differential graded algebra given by
\begin{equation}\label{koszul}
K =  K(r_1) \otimes \dotsb \otimes K(r_n),
\end{equation}
according to a prescription similar to \eqref{monoid}, but for
chain complexes. As a graded algebra, $K$ is the exterior algebra
$\Lambda(e_1,\ldots,e_n)$. Together with the augmentation
$\epsilon= \otimes_{i=1}^n \epsilon_{r_i}$,
$K$ is an $R$--free resolution of $R/I$ of length $n$
(\cite{M}*{Theorem 16.5}).

Changing to two dimensions, we start by realizing $K(r)$ and
$R[x]$, for a non-zero divisor $r\in R$ and a free variable $x$,
as multiplicative double complexes, by assigning $e$ and $x$
bidegrees $(0,1)$ and $(-1,1)$ respectively. So $K(r)$ is
concentrated on the $y$--axis and $R[x]$ on the secondary
diagonal. Of course, all the differentials of $R[x]$ are trivial.

The multiplicative double complex $\KK(r)$ corresponding to $K(r)$
is now given as follows. We take $K(r)\otimes R[x]$ and define
horizontal differentials $d^h$ by setting $d^h_{-k,k+1}(e\otimes
x^k) = 1 \otimes x^{k+1}$. The other components of $d^h$ are
necessarily trivial. To see that $d^h$ is compatible with the
multiplicative structure canonically defined on $K(r)\otimes
R[x]$, we have to check equation \eqref{diffbigradedh}. We do the
calculation for $c=e\otimes x^k$ and $d=e\otimes x^l$. The
left-hand side of the equation is trivial, and the right-hand side
is given by
\begin{multline*}
   \mu\bigl((1 \otimes x^{k+1})\otimes(e\otimes x^l)\bigr) + (-1)^{-k}\,
\mu\bigl((e\otimes x^k)\otimes (1\otimes x^{l+1})\bigr) = \\
(-1)^{k+1}\, e\otimes x^{k+l+1} + (-1)^{-k}\, e\otimes x^{k+l+1} =
0.
\end{multline*}
The vertical differentials are induced by the ones of $K(r)$.
Explicitly, its non-trivial components are $d^v_{-k, k+1}(e\otimes
x^k) = r\otimes x^k$.

\begin{definition}
The {\em extended Koszul complex}\/ $\KK$ associated to the
regular sequence $r_1, \ldots, r_n$ is the multiplicative double
complex
\[ \KK = \KK(r_1)\otimes \dotsb \otimes \KK(r_n).\]
\end{definition}

Let us describe the bigraded algebra $U(\KK)$ underlying $\KK$.
Underlying a building block $\KK(r_i)$ is the bigraded algebra
\[ U(\KK(r_i; x_i)) = \Lambda(e_i)\twist R[x_i].\]
Note that the elements $e_i$ and $x_i$ {\em anti}commute (we have
used this in the verification above). Moreover, we have the
identifications (of bigraded algebras)
\begin{align*}
\Lambda(e_1)\twist\Lambda(e_2) & = \Lambda(e_1, e_2)\\
R[x_1]\twist R[x_2] & = R[x_1, x_2].
\end{align*}
Consequently, $U(\KK)$ is isomorphic to a twisted tensor product
of an exterior and a polynomial algebra, concentrated on the
$y$--axis and the secondary diagonal respectively, by means of a
composition of symmetry isomorphisms. For $n=2$, we have for
instance
\begin{align*}
U(\KK) & = U(\KK(r_1))\twist U(\KK(r_2)) = \Lambda(e_1) \twist
R[x_1]\twist \Lambda(e_2)\twist R[x_2] \\ & \cong
\Lambda(e_1)\twist\Lambda(e_2)\twist R[x_1]\twist R[x_2] =
\Lambda(e_1,e_2)\twist R[x_1, x_2],
\end{align*}
given on elements by
\begin{equation*}
e_1^{j_1}\otimes x_1^{i_1}\otimes e_2^{j_2}\otimes x_2^{i_2}
\mapsto (-1)^{j_2 i_1} e_1^{j_1}\wedge e_2^{j_2}\otimes x_1^{i_1}
x_2^{i_2},
\end{equation*}
for $j_1, j_2\in\{0,1\}$ and $i_1, i_2\geq 0$. For arbitrary $n$,
an isomorphism
\begin{equation}\label{underlyingalgebra}
U(\KK) \cong \Lambda(e_1, \ldots, e_n)\twist R[x_1,\ldots, x_n]
\end{equation}
is given by
\begin{equation}
e_1^{j_1}\otimes x_1^{i_1}\otimes\dotsb \otimes e_n^{j_n}\otimes
x_n^{i_n} \mapsto (-1)^{(i,j)} e_1^{j_1}\wedge \dotsb \wedge
e_n^{j_n}\otimes x_1^{i_1}\dotsm x_n^{i_n}
\end{equation}
for $ j_1, \ldots, j_n\in\{0, 1\}$, $i_1,\ldots,i_n\geq 0$, where
$(i,j)$ is defined as \[ (i,j) = j_2 i_1 + (j_3 i_2 + j_3 i_1) +
\dotsb + (j_n i_{n-1} + \dotsb + j_n i_1).\] Under this
isomorphism, the diagonal line $p+q=k$ of $\KK$ corresponds to the
(right) $R[x_1, \ldots, x_n]$--submodule of $\Lambda(e_1,\ldots,
e_n) \twist R[x_1, \ldots, x_n]$ generated by the homogeneous
elements of degree $k$ of $\Lambda(e_1, \ldots, e_n)$, for
instance. We will refer to elements of $\KK$ by means of
\eqref{underlyingalgebra}. The element $(e_1\otimes 1)\otimes
(1\otimes x_2^4) \otimes (e_3\otimes 1)$ for example is written as
$- e_1\wedge e_3\otimes x_2^4$.

The column-wise filtration of $\KK$ arises naturally as a
filtration by ideals. Namely, define $F^s(\KK)$ for $s\geq 0$ to
be the (left) ideal generated by $J^s$, where $J$ is the maximal
ideal $J=(x_1,\ldots, x_n) \lhd R[x_1, \ldots, x_n]$, so
\[ F^s(\KK)  = \KK \. J^s.\]

We denote the quotients of $\KK$ by these ideals by $\KK/s$ and
the components of the associated graded by $\QQ^s$,
\[ \KK/s = \KK / F^s(\KK), \quad \QQ^s= F^{s}(\KK) /F^{s+1}(\KK).\]
Put more directly, $\KK/s$ consists of the first $s$ columns (with
the differentials from $\KK$) and $\QQ^s$ of the $s$th column of
$\KK$; however, $\QQ^s$ is still a double complex. Its
condensation agrees with the $s$th column of $\KK$ up to a shift,
namely $(\QQ^s)^\bullet = \Sigma^{-s}\KK^s$, where $\Sigma$
denotes the suspension functor, defined on a chain complex $C$ as
$(\Sigma C)_k = C_{k-1}$, with differentials $(\Sigma d)_k =
-d_k$.

Let us determine the homology of the columns and the rows of
$\KK$. The columns can be identified as direct sums of the Koszul
complex $K(r_1, \ldots, r_n)$, indexed by the homogeneous
monomials in $x_1, \ldots, x_n$ of degree $s$. Namely, we have
\begin{equation}\label{columnkoszul}
\KK^s \cong K\otimes J^s/J^{s+1},
\end{equation}
essentially by definition:
\begin{align*}
 \KK^s \ & = \bigoplus_{i_1+\dotsb+ i_n=s} \KK(r_1)^{i_1}
 \otimes\dotsb\otimes \KK(r_n)^{i_n}\\ &=
\bigoplus_{i_1+\dotsb+ i_n=s} K(r_1)\otimes
Rx_1^{i_1}\otimes\dotsb\otimes K(r_n)\otimes R x_n^{i_n}\\ &\cong
\bigoplus_{i_1+\dotsb+ i_n=s} K \otimes R
x_1^{i_1}\otimes\dotsb\otimes R x_n^{i_n} \\ & \cong \quad
K\otimes J^s/J^{s+1}.
\end{align*}
This implies
\begin{equation}
\label{columnhomology} H_p(\KK^s)\cong
\begin{cases}J^s/J^{s+1} & \text{if }p=s,\\0 &\text{otherwise}.\end{cases}
\end{equation}
The rows of $\KK(r_i)$ are all chain complexes of free
$R$--modules; other than the zeroth one, which consists of $R$
(concentrated in degree zero), they are exact. Therefore, using
property \eqref{rows}, the K\"unneth theorem implies that
\begin{equation}
\label{rowhomology} H_p(\KK_t) =
\begin{cases} R &\text{if } p=t=0,\\ 0 &\text{otherwise}.\end{cases}
\end{equation}

Next, we define an augmentation on $\KK^\bullet$. Note that
$d_1\:\KK^\bullet_1\to \KK^\bullet_0$ is up to sign given by
$d_1(e_j\otimes f)=\pm (r_j-x_j) f$, for $f\in R[x_1,\ldots,
x_n]$. Hence the evaluation map
\[\epsilon\: R[x_1, \ldots, x_n] \lra R, \quad \epsilon(x_i)= r_i\]
defines an augmentation. As $\epsilon$ is compatible with the
filtrations given by powers of the ideals $J\lhd R[x_1,\ldots,
x_n]$ and $I\lhd R$ respectively, it induces augmentations
\begin{alignat*}{2} \epsilon^s\: (F^s(\KK))^\bullet_0 \lra I^s,
\quad \epsilon_s\: (\KK/s)^\bullet_0 \lra R/I^s,
\quad\epsilon^s_{s+1}\: (\QQ^s)^\bullet_0 \lra I^s/I^{s+1}
\end{alignat*} for the complexes
$(F^s(\KK))^\bullet$, $(\QQ^s)^\bullet$, $(\KK/s)^\bullet$
respectively.

\begin{proposition}[{compare \cite{B}*{Theorem 1.3}}]
\label{resolutions} For $s\geq 1$, the differential graded
algebras $(F^s(\KK))^\bullet$, $(\QQ^s)^\bullet =
\Sigma^{-s}(K\otimes J^s/J^{s+1})$ and $(\KK/s)^\bullet$ provide
$R$--free resolutions
\begin{align*}
& (F^s (\KK))^\bullet \xrightarrow{\epsilon^s} I^s \lra 0, \\
& (\KK/s)^\bullet \xrightarrow{\epsilon_s} R/I^s \lra 0\\
& (\QQ^s)^\bullet \xrightarrow{\epsilon^s_{s+1}} I^s/I^{s+1} \lra
0
\end{align*}
of $I^s$, $R/I^s$ and $I^s/I^{s+1}$ respectively.
\end{proposition}

\begin{proof}
Filtering the respective double complexes by rows gives rise to
Cartan--Eilenberg spectral sequences (\cite{CE}*{XV, \S 6}). They
converge because the filtrations are bounded below and exhaustive,
as the components of the condensation are defined as direct sums
(see \cite{We}*{5.6}). The $E^1$--term, given by the homology of
the columns, is concentrated on the diagonal $p+q=0$, as a
consequence of \eqref{columnhomology}.

Therefore, it only remains to check exactness of the complexes in
degree zero. This is clear; the augmentations induce isomorphisms
\begin{align*}
& H_0((F^s \KK)^\bullet) = J^s/(x_1-r_1, \ldots, x_n-r_n)
\xrightarrow{\cong} I^s,\\
& H_0((\KK/s)^\bullet) = R[x_1,\ldots, x_n]/(J^s, x_1-r_1, \ldots
x_n-r_n) \xrightarrow{\cong} R/I^s,\\
& H_0((\QQ^s)^\bullet) = J^s/(J^{s+1}, x_1-r_1, \ldots, x_n-r_n)
\xrightarrow{\cong} I^s/I^{s+1}. \qedhere
\end{align*}
\end{proof}

\begin{remark}\label{generalresolution}
A similar argument shows that there are $R$--free resolutions of
any subquotient of $R$ of the form $I^s/I^t$, for $t>s\geq 0$,
given by
\[ (F^s(\KK)/F^{t}(\KK))^\bullet \xrightarrow{\epsilon^s_{t}}
 I^s/I^{t} \lra 0, \]
where $\epsilon^s_{t}$ is induced by $\epsilon$.
\end{remark}

\begin{remark}\label{cover}
Note that the canonical inclusion $(\QQ^s)^\bullet \to
(\KK/s)^\bullet$ covers (via the augmentations) the inclusion
$I^s/I^{s+1}\to R/I^s$; similarly, the canonical projection $(F^s
(\KK))^\bullet \to (\QQ^s)^\bullet$ covers the projection $I^s\to
I^s/I^{s+1}$.
\end{remark}

\section{Computation of Tor--groups}

Tensoring the Koszul resolution $K$ with $R/I$ kills all the
differentials. Therefore Proposition \ref{resolutions} implies
that
\begin{equation}\label{recognition}
 \Tor_*(R/I, I^s/I^{s+1}) = R/I\otimes \Sigma^{-s} (\KK^s\otimes
J^s/J^{s+1})
\end{equation}
(recall that $J$ denotes the ideal $J=(x_1, \ldots, x_n)\lhd
R[x_1, \ldots, x_n]$). Together with \eqref{columnkoszul} we find
the following well-known result:

\begin{proposition}\label{homologyk}
For $s\geq 0$, we have
\[ \Tor_*(R/I, I^s/I^{s+1}) =  \Lambda_{R/I}(e_1, \ldots, e_n)\otimes
J^s/J^{s+1}.\] For $s=0$, this is an identity of
algebras.\hfill\qed
\end{proposition}

\begin{remark}
The statement for $s>0$ follows in fact directly from the case
$s=0$, as a consequence of the splitting (of $R$--modules)
\[ I^s/I^{s+1} \cong R/I\otimes J^s/J^{s+1} \cong \bigoplus_{x\in
V_s} R/I\ x,\] where $V_s$ is the set of monomials of degree $s$
in $x_1, \ldots, x_n$ (\cite{M}*{Theorem 16.2}).
\end{remark}

It is convenient to introduce a  reduced version of $\Tor$ for
$R/I^s$. For this, let us first note that the projection $R\to
R/I^s$ induces a split monomorphism \[R/I\cong \Tor_*(R/I, R)\lra
\Tor_*(R/I, R/I^s).\] This can be seen as follows. Replace $R/I$
and $R/I^s$ by projective resolutions $P$ and $Q$ which are
differential graded algebras, e.g.\@ $P=(\KK/1)^\bullet$ and $Q =
(\KK/s)^\bullet$. Then the map $\Tor_*(R/I, R)\to \Tor_*(R/I,
R/I^s)$ is induced by $P \to P\otimes Q$ and split by $P\otimes Q
\to P\otimes P \to P$, where the first map is the identity on $P$
tensored with a lift of the projection $R/I^s\to R/I$. Defining
\[ \widetilde \Tor_*(R/I, R/I^s) = \coker(R/I \lra \Tor_*(R/I, R/I^s)),\]
we therefore have
\[ \Tor_*(R/I, R/I^s) \cong R/I \oplus \widetilde \Tor_*(R/I, R/I^s). \]
We need a few notations. For $s\geq 0$, let $\partial_s$ be the
connecting homomorphism \[\partial_s:\Tor_{*+1}(R/I, I^s/I^{s+1})
\lra \Tor_*(R/I, I^{s+1}/I^{s+2})\] associated to the short exact
sequence
\begin{equation}\label{sesconn}
0 \lra I^{s+1}/I^{s+2}\lra I^s/I^{s+2} \lra I^s/I^{s+1}\lra 0,
\end{equation}
$i_s$ the inclusion $I^{s}/I^{s+1} \to R/I^{s+1}$ and $p_s$ the
projection $I^s\to I^s/I^{s+1}$. We will omit the index $s$ if it
is clear from the context. The projection $\Tor_*(R/I, R/I^s) \to
\widetilde \Tor_*(R/I, R/I^s)$ is denoted by $\pi$.

\begin{theorem}[{\cite{B}*{Lemma 2.1, Proposition 2.2}}]
\label{thmquotients} For $s>1$, the sequence of graded
$R/I$--modules
\[ \Tor_{*+1}(R/I, I^{s-2}/I^{s-1}) \xrightarrow{\partial}
\Tor_{*}(R/I, I^{s-1}/I^s) \xrightarrow{\pi i_*} \widetilde
\Tor_*(R/I, R/I^s) \lra 0\] is exact; moreover, $\widetilde
\Tor_*(R/I, R/I^s)$ is free over $R/I$. The product structure on
$\widetilde \Tor_*(R/I, R/I^s)$, induced by the $R$--algebra
structure on $R/I^s$, is trivial.
\end{theorem}

To simplify notation, we abbreviate $\Tor_*(R/I,-)$ by $H_*(-)$ in
the following and refer to it as homology; similarly, $\widetilde
H_*(-)$ stands for $\widetilde\Tor_*(-)$.

\begin{proof}
Making use of the free resolution $(\KK/s)^\bullet$ of $R/I^s$
constructed in the previous section, we can compute $H_*(R/I^s)$
as $H_*(R/I\otimes(\KK/s)^\bullet$). The filtration defined for
$\KK$ induces one on $\KK/s$, of the form
\begin{equation}
\label{filtrationdouble} 0 = F^s\subset F^{s-1} \subset \dotsb
\subset F^0 = \KK/s.
\end{equation} Setting $G^s= R/I\otimes (F^s)^\bullet$, this in turn
gives rise to the filtration
\begin{equation}\label{filtrationcondensation}
0=G^s\subset G^{s-1} \subset \dotsb  \subset G^0 =
R/I\otimes(\KK/s)^\bullet,
\end{equation}
which is a filtration of the differential graded algebra
$R/I\otimes(\KK/s)^\bullet$ by ideals. It determines a
Cartan--Eilenberg spectral sequence converging to $H_*(R/I^s)$, by
the same argument as in the proof of Proposition
\ref{resolutions}. The $E^1$--term is given by the homology of the
columns, which we have identified --- up to a shift --- as free
resolutions of the quotients $I^p/I^{p+1}$, so $(E^1)^p_{*+p} =
H_*(I^p/I^{p+1})$. The first differential is given by the
horizontal differential, therefore \eqref{rowhomology} implies
\[
(E^2)_*^p  =
\begin{cases}
\coker((d^h)^{s-2}\:(E^1)_*^{s-2}\lra (E^1)_*^{s-1}) & \text{if } p= s-1,\\
R/I &\text{if } p=0,\\
0 & \text{otherwise}.
\end{cases}
\]
For dimensional reasons, the spectral sequence collapses at $E^2$.
Consequently, the composition
\[(E^2)^{s-1}_{*+s-1}=(E^\infty)^{s-1}_{*+s-1} = H_*(G^{s-1})
\hookrightarrow H_{*}(R/I^s) \xrightarrow{\pi} \widetilde
H_*(R/I^s)\] of the edge homomorphism with the projection $\pi$ is
an isomorphism, so that we have an exact sequence
\[ (E^1)_{*+s-1}^{s-2}\lra (E^1)_{*+s-1}^{s-1} \lra \widetilde H_*(R/I^s) \lra 0.\]
Remark \ref{cover} shows that the map $(E^1)_{*+s-1}^{s-1} \to
\widetilde H_*(R/I^s)$ can be identified with $\pi
i_*\:H_{*}(I^{s-1}/I^s)\to \widetilde H_*(R/I^s)$.

It remains to identify $R/I\otimes d^h\:(E^1)_*^{s-2}\to
(E^1)_*^{s-1}$ as a connecting homomorphism. The free resolutions
of the terms in the short exact sequence \[ 0 \lra I^{s-1}/I^s\lra
I^{s-2}/I^s \lra I^{s-2}/I^{s-1} \lra 0,\] described in
Proposition \ref{resolutions} and Remark \ref{generalresolution},
fit into a short exact sequence of chain complexes
\[ 0 \lra (\QQ^{s-1})^\bullet \lra (F^{s-2}(\KK)/F^s(\KK))^\bullet \lra
(\QQ^{s-2})^\bullet \lra 0. \] Going through the definition of the
connecting homomorphism, we find that $\partial_{s-2}$ is given by
$R/I \otimes (d^h)^{s-2}$.

For the second statement, it suffices to observe that the image of
$\partial$ is a free submodule of the free $R/I$--module
$H_*(I^{s-1}/I^s)$. This is true because $\partial=R/I\otimes d^h$
maps basis elements of $R/I\otimes\KK$ to sums of basis elements.

For the determination of the multiplicative structure of
$H_*(R/I^s)$, recall from Section \ref{multiplicative} that the
product on $\KK/s$ induces one on $(\KK/s)^\bullet$ and hence on
$H_*(R/I^s)=H_*(R/I\otimes(\KK/s)^\bullet)$. The latter one indeed
{\em is} the canonical internal product defined on $H_*(R/I^s)$
(see \cite{M2}*{Corollary VIII.2.3}). As the product on $\KK/s$ is
compatible with the filtration \eqref{filtrationdouble}, so are
the induced ones on $R/I\otimes(\KK/s)^\bullet$ and on
$H_*(R/I^s)$.
Now we have seen above that $\widetilde H_*(R/I^s)$ is
concentrated in $H_*(G_{s-1})$, on which the product is trivial.
\end{proof}

We can express $H_*(I^s)$ in a similar manner (the proof is
completely analogous):

\begin{theorem}\label{thmideals}
For $s\geq 0$, there is an exact sequence
\[ 0 \lra \Tor_*(R/I, I^s) \xrightarrow{p_*} \Tor_*(R/I, I^s/I^{s+1})
\xrightarrow{\partial} \Tor_{*-1}(R/I, I^{s+1}/I^{s+2});\]
\item\label{elements}
moreover, $\Tor_*(R/I, I^s)$ is free over $R/I$. The product
structure on $\Tor_*(R/I, I^s)$, induced by the $R$--algebra
structure on $I^s$, is trivial.\hfill\qed
\end{theorem}

\begin{remark}\label{tor0}
In particular, we have $R/I\otimes I^s \cong I^s/I^{s+1}$.
\end{remark}

It is quite easy to make the connecting homomorphism
\[\partial\:H_{*+1}(I^{s}/I^{s+1}) \to H_{*}(I^{s+1}/I^{s+2})\]
explicit. In the following, we slightly abuse notation and denote
basis elements of
\[
\bigoplus_{s\geq 0} H_*(I^s/I^{s+1}) =
\bigoplus_{s\geq 0}\, \Sigma^{-s}(R/I\otimes\KK)^s
\]
(see \eqref{recognition}) as expressions in the variables $e_i$
and $x_i$, which are in fact generators of the components of
$\KK$. By means of illustration, consider first $\partial\:
H_{*+1}(R/I) \to H_*(I/I^2)$. The unit is mapped to zero and the
elements $e_i$ to $x_i$. For the elements of higher degree, we
have
\begin{align*}
 \partial(e_1 \wedge e_2) & = - e_2\otimes x_1 + e_1\otimes x_2\\
\partial(e_1\wedge e_2\wedge e_3) & = e_2\wedge e_3\otimes x_1 - e_1\wedge
e_3\otimes x_2 + e_1\wedge e_2\otimes x_3,
\end{align*}
and so on. Considering $R/I\otimes\KK$ under the isomorphism
\eqref{underlyingalgebra} as a right module over $R/I[x_1, \ldots,
x_n]$, $\partial$ is a linear map. So we have for instance \[
\partial(e_1\wedge e_2 \otimes x_1) = -e_2 \otimes x_1^2  + e_1 \otimes x_1 x_2.\]
In the general case, we have the following formula:

\begin{proposition}\label{connectinghom}
The map $\partial\:H_{*+1}(I^{s}/I^{s+1}) \to
H_{*}(I^{s+1}/I^{s+2})$ is given by
\[\partial( e_{i_1}\wedge\dotsb\wedge e_{i_l}\otimes f) =
\sum_{j=1}^{l} (-1)^{j+l} \ e_{i_1}\wedge\dotsb\wedge \widehat
e_{i_j}\wedge\dotsb\wedge e_{i_{l}} \otimes x_{i_j} f,\] where
$\{i_1, \ldots, i_l\}\subset\{1, \ldots, n \}$, $f$ is a monomial
in $x_1,\ldots, x_n$ of degree $s$, and the hat indicates that the
entry underneath should be omitted.\hfill\qed
\end{proposition}

\begin{remark}\label{explicitbasis}
Identifying $H_*(I^s)$ with its image under $p_*$ and recalling
that $\ker\partial_s=\im\partial_{s-1}$ for $s\geq 1$, the
proposition gives an explicit description of a basis of
$H_*(I^s)$.
\end{remark}

\begin{remark} Because all the
$\Tor$--groups we have computed are free $R/I$--modules, we have
the $\Ext$--groups for free, because of the following fact. If $A$
is an $R$--module such that $\Tor_*^R(R/I, A)$ is free over $R/I$,
there is a duality isomorphism
\[ \Ext^*_R(A, R/I) \cong \Hom^*_{R/I}(\Tor_*^R(R/I, A), R/I).\] It
arises as an edge homomorphism of a Cartan--Eilenberg spectral
sequence, which collapses under this condition (see
\cite{CE}*{XVI, \S 6, Case 3}).
\end{remark}

\section{An exterior multiplication}
As mentioned in the introduction, we aim to study
$\bigoplus_{s\geq 0} \Tor_*(R/I,I^s)$ as a bigraded
$R/I$--algebra. We abbreviate $\Tor_*(R/I,-)$ by $H_*(-)$, as
before.

Ordinary ring multiplication induces pairings $I^s\otimes I^t\to
I^{s+t}$. Taken all together, these give rise to a graded
$R$--algebra structure on
\[ B_I^*(R) = R \oplus I\oplus I^2 \oplus \dotsb\]
(called the blowup algebra of $I$ in $R$ in \cite{E}*{5.2}). On
the other hand, we have the graded ring associated to the
$I$--adic filtration
\[ \gr_I^*(R) = R/I \oplus I/I^2 \oplus I^2/I^3 \oplus\dotsb \]
These product structures induce exterior multiplications on
homology, giving each $H_*(B_I^*(R))$ and $H_*(\gr_I^*(R))$ the
structure of a bigraded $R/I$--algebra. We may compute these using
the multiplicative double complexes deduced from $\KK$.
Proposition \ref{homologyk} immediately implies
\begin{proposition}
There is an isomorphism of bigraded algebras
\[ \Tor_*(R/I, \gr_I^*(R))\cong \Lambda_{R/I}(e_1,
\dots, e_n)\twist R/I[x_1,\ldots, x_n].\hfill\qed\]
\end{proposition}
\noindent The projection $I^s\to I^s/I^{s+1}$ induce a map of
graded $R$--algebras
\[ p\: B^*_I(R) \lra \gr_I^*(R).\]
On $\Tor$, this induces by Theorem \ref{thmideals} a monomorphism
of bigraded $R/I$--algebras
\begin{equation}
p_*\: H_*(B_I^*(R)) \lra  H_*(\gr_I^*(R));
\end{equation}
we identify $H_*(B_I^*(R))$ with its image under $p_*$ in the
following. By Remark \nolinebreak \ref{tor0}, \[\Tor_0(R/I,
B_I^*(R))\cong R/I[x_1, \ldots, x_n],\] hence we can consider
$H_*(B_I^*(R))$ as a (bigraded) algebra over $R/I[x_1, \dots,
x_n]$. Note that $R/I[x_1,\ldots, x_n]$ is {\em not} contained in
the centre of $H_*(B_I^*(R))$.

\begin{proposition}
Over $R/I[x_1,\ldots, x_n]$, the bigraded algebra $\Tor_*(R/I,
B_I^*(R))$ is generated by the basis elements
\[ a_{(i_0, \ldots, i_k)} = \partial(e_{i_0}\wedge\dotsb\wedge e_{i_k})\]
of $\Tor_k(R/I, I)$ for $0<k<n$, where $\{i_0, \ldots, i_k\}$ runs
through the subsets of $\{1,\ldots, n\}$ of cardinality $k+1$ and
we assume that $1\leq i_0< \dotsb < i_k \leq n$.
\end{proposition}

\begin{proof}
It is clear from the short exact sequence  $0 \to I \to R \to R/I
\to 0$ and Proposition \ref{homologyk} that the $a_{(i_0, \ldots,
i_k)}$ defined in the statement generate $H_*(I)$ as an
$R/I$--module. Now we claim that a set of generators of the $s$th
column $H_*(I^s)$ of $H_*(B_I^*(R))$ is given by multiplying all
monomials of degree $s$ in $x_1, \ldots, x_n$ with all the
$a_{(i_0, \ldots, i_k)}$. Namely, we know that, for $s>0$,
\begin{align*}
H_*(I^s) & \cong \ker(\partial_{s}\: H_{*}(I^s/I^{s+1}) \lra
H_{*-1}(I^{s+1}/I^{s+2}))\\
& = \im(\partial_{s-1}\: H_{*+1}(I^{s-1}/I^s) \lra
H_*(I^s/I^{s+1})),
\end{align*}
(the first equality is Theorem \ref{thmideals}, the second is
equation \eqref{rowhomology} together with the recognition of
$R/I\otimes d^h$ as the connecting homomorphism $\partial$, in the
proof of Theorem \ref{thmquotients}), and Proposition
\ref{connectinghom} implies that $\im\partial_{s-1} =
J^{s-1}/J^s\. \im\partial_0$.
\end{proof}

We finish by giving some examples for small $n$; we abbreviate
$H_*(B_I^*(R))$ by $A$, $R/I$ by $E$ and $R/I[x_1,\ldots, x_n]$ by
$P$.
\begin{itemize}
\item $n=1$: Clearly, we have $A\cong P = E[x]$.
\item $n=2$: The basis element $a_{12}= - e_2 x_1+e_1 x_2$ of $H_1(I)$
generates a free copy of $P$. More precisely,
$A\cong\Lambda_E(a_{12})\twist P$, as bigraded algebras.
\item $n=3$: Among the basis elements $a_{12}$, $a_{13}$ and $a_{23}$
we have the relation
\begin{equation*}\label{prelation}
x_1\. a_{23} + x_2\. a_{13} + x_3\. a_{12} = 0
\end{equation*}
over $P$. As an example of a product, we have $a_{12}\. a_{23} = -
x_2\cdot a_{123}$.
\item $n=4$: In addition to four relation of the type above, there is also
\[ x_1\. a_{234} + x_2\. a_{134} + x_3\. a_{124} + x_4\. a_{123}
=0.\] As a product, we have $a_{123}\. a_{234} = x_2 x_3\.
a_{1234}$.
\end{itemize}

\end{document}